\documentclass[11pt]{article}
\usepackage[T1, T2A]{fontenc}
\usepackage[utf8]{inputenc}
\usepackage[russian,english]{babel}
\usepackage{amssymb,amsfonts,amsmath,mathtext,mathrsfs}
\usepackage[matrix,arrow,curve]{xy}

\def\cF{{\cal F}}
\def\cG{{\cal G}}
\def\cB{{\cal B}}

\def\cH{{\cal H}}
\def\cP{{\cal P}}
\def\CC{{{\mathbb{C}}}}
\def\FF{{\mathbb F}}

\usepackage{color}

\definecolor{MyGrayc}{gray}{0.75}
\definecolor{MyGreen}{rgb}{0.0,0.39368872549019607,0.0}
\definecolor{MyRed}{rgb}{0.8, 0, 0}
\definecolor{MyBlue}{rgb}{0,0,0.8}
\definecolor{my-color-chocolate}{rgb}{0.8,0.5,0.4}
\definecolor{my-color-dark-magenta}{rgb}{0.5472273284313726,0.0,0.5472273284313726}

\usepackage[unicode=true,colorlinks=true,linkcolor=blue,filecolor=magenta,urlcolor=cyan,
bookmarksopen=true,
,citecolor=my-color-chocolate]{hyperref}

\usepackage[a4paper,hmargin=3cm,vmargin=2.5cm]{geometry}
\usepackage{theorem}
\theoremstyle{break}

\newtheorem{Th}{Theorem}[section]
\newtheorem{Prop}{Proposition}[section]
\newtheorem{Corollary}[Prop]{Corollary}
\newtheorem{Le}[Prop]{Lemma}

\newtheorem{Ax}[Prop]{Axiom}
{\theorembodyfont{\rmfamily}
  \newtheorem{Def}[Prop]{Definition}
  \newtheorem{Ex}[Prop]{Example}
  \newtheorem{Rem}[Prop]{Remark}
  
}

\newenvironment{Proof}{\par\noindent{\bf Proof.\\}}
{\hfill$\scriptstyle\blacksquare$}

\usepackage{ifthen}
\def\enableComments{0}
\newcommand{\commentWhenEnable}[1]%
{\ifthenelse{\equal{\enableComments}{1}}%
  {#1}{}}

\newcounter{marpCounter}

\usepackage{cancel}

\newcommand{\comment}[1]{}

\def\enableTODOlist{0}
\newcommand{\TODO}[1]{\ifthenelse{\equal{\enableTODOlist}{1}}%
  {#1}%
  {}}

\usepackage{epsf}
\usepackage{graphicx}

\title{An extension of Stanley's chromatic symmetric function to binary delta-matroids}
\author{M.~Dudina\thanks{Skolkovo Institute of Science and Technology}, V.~Zhukov\thanks{National Research University Higher School of Economics}}

\usepackage{color}\definecolor{RED}{rgb}{1,0,0}\definecolor{myblue}{rgb}{0,0.5,1}

\begin{document}
\maketitle
\section{Introduction}
The chromatic polynomial is a well-known and extensively studied graph invariant
with values in a polynomial ring in one variable~$t$. It enumerates proper colorings
of the vertices of a simple graph in~$t$ colors.
Stanley's symmetrized chromatic polynomial is a generalization of the ordinary chromatic
polynomial to a graph invariant with values in a ring of polynomials in infinitely many
variables. The ordinary chromatic polynomial is a specialization of the Stanley's one.

Our goal is to extend Stanley's chromatic polynomial to embedded graphs.
In contrast to well-known extensions of, say, the Tutte polynomial
from abstract to embedded graphs~\cite{chun2014matroids}
, we do not treat an embedded graph as
an abstract graph endowed with additional information about the embedding.
Instead, we consider the binary delta-matroid associated to an embedded graph
and define the extended Stanley's chromatic polynomial as an invariant
of delta-matroids.

We show that, similarly to Stanley's symmetrized chromatic polynomial of graphs,
which satisfies $4$-term relations for simple graphs and determines in this way
a knot invariant,
extended Stanley's chromatic polynomial of binary delta-matroids
we define satisfies the $4$-term relations for binary delta-matroids~\cite{LZ}
and determines, therefore, an invariant of {\em links}.\\
In Sec.~\ref{major-notions}  the necessary background on graphs and delta-matroids is provided, together with the structures of the related Hopf algebras. The concept of combinatorial Hopf algebras is reviewed in Sec.~\ref{CHA}.
 We start Sec.~\ref{Chr} with giving two major definitions of Stanley's  symmetric chromatic function for graphs.
 Then, we extend the second definition, the one using characters in Hopf algebras,
 to define a symmetric chromatic function for delta-matroids. In Sec.~\ref{4T},
 the newly defined chromatic symmetric function is proven to satisfy the extended
 four-term relation on binary delta-matroids. In Sec.~\ref{Prim} we compute the values of the new invariant on the space of primitive elements corresponding to even binary delta-matroids.

\section{Major notions}

Here we briefly describe the structure of underlying Hopf algebras for graphs, framed graphs, delta-matroids and binary delta-matroids.
\subsection{Graphs and Delta-Matroids}
\label{major-notions}
\subsubsection{Hopf algebras of graphs and framed graphs}

The Hopf algebra~$\cG$ of graphs here is the graded vector space
$$
\cG=\cG_0\oplus\cG_1\oplus\dots,
$$
with~$\cG_n$ being the vector space spanned over~$\CC$ by isomorphism
classes of graphs having~$n$ vertices.
The multiplication in~$\cG$ is induced by the disjoint union of graphs,
and the comultiplication $\mu:\cG\to\cG\otimes\cG$ acts on a graph~$G$ as
$$
\mu: G\to \sum_{{U\sqcup W=V(G)}}G|_U\otimes G|_W,
$$
where the summation runs over all partitions of the set of vertices~$V(G)$
into a disjoint union of two subsets, $G|_U$ denoting the subgraph of~$G$
induced by the subset~$U\subset V(G)$. (The comultiplication is often denoted by~$\Delta$,
but in the setting of the present paper the symbol~$\Delta$ will usually denote the symmetric difference
of sets.)

In addition to the Hopf algebra of graphs~$\cG$ we will consider three other Hopf algebras:
\begin{itemize}
\item the Hopf algebra~$\cG^f$ of framed graphs;

\item the Hopf algebra~$\cB$ of binary delta-matroids;

\item the Hopf algebra~$\cB^e$ of even binary delta-matroids.
\end{itemize}
Note that the relationship between graphs and framed graphs is the same
as the one between even binary delta-matroids and delta-matroids.

A {\em framed graph\/} is a simple graph~$G$ endowed with a {\em framing},
which is a mapping $V(G)\to\{0,1\}$. Similarly to ordinary graphs,
a framed graph can be represented by its {\em adjacency matrix}.
The columns and the rows of the adjacency matrix are numbered
by the vertices of the graph.
The non-diagonal entries of the adjacency matrix are equal to~$1$
provided the corresponding vertices are connected by an edge and~$0$ otherwise.
The diagonal entries are equal to the framing of the corresponding vertex.

The {\em Hopf algebra of framed graphs\/} (introduced in~\cite{L06}) is
$$
\cG^f=\cG^f_0\oplus\cG^f_1\oplus\cG^f_2\oplus\dots,
$$
where $\cG^f_n$, $n=0,1,2,\dots$, is the vector space spanned by
isomorphism classes of framed graphs with~$n$ vertices.
Note that the grading~$1$ vector space~$\cG^f_1$ is spanned by the
two framed graphs having a single vertex, one of them with the framing~$0$,
the other one with the framing~$1$.
The multiplication $\cG^f\otimes\cG^f\to\cG^f$ and the comultiplication
$\cG^f\to\cG^f\otimes\cG^f$ are defined in the same vein as for the
Hopf algebra~$\cG$.
\subsubsection{Hopf algebra of binary delta-matroids}
\label{sec:--simple-grp-and-even-binary-d-matroids}

Delta-matroids were introduced by A.~Bouchet~\cite{bouchet1988representability}. Below, we mainly follow the
approach and terminology from~\cite{chun2014matroids}.

\begin{Def}
  A {\em set system} is a pair $(E; \Phi)$, where $E$ is an arbitrary finite set,
  and  $\Phi\subset 2^E$ is a set of subsets of $E$.

  The set $E$ is called the  {\em ground set} of the set system $(E; \Phi)$,
  and the  elements of  $\Phi$ are called  the {\em feasible sets } of the set system $(E; \Phi)$.

  A set system  $(E; \Phi)$ is said to be
  {\em proper} provided $\Phi$ is non-empty. (Below, we consider proper systems only.)

Two set systems $(E;\Phi)$ and $(E';\Phi')$ are said to be \emph{isomorphic}
if there is a one-to-one mapping $E\to E'$ transforming~$\Phi$ to~$\Phi'$.

\end{Def}

\begin{Def}
  A proper set system $(E; \Phi)$ is a {\em delta-matroid} if the following axiom {\rm(}the {\em
  Symmetric Exchange Axiom, SEA}{\rm)} is satisfied:
  \begin{Ax}[SEA]
  \label{SEA}
    For any two feasible sets $X,Y\in\Phi$  and for any element $a\in X\Delta Y$
    there exists an element $b\in X\Delta Y$ {\rm(}which is allowed to be equal to $a${\rm)}
    such that $X\Delta\{a, b\}\in\Phi$ {\rm(}in the case $a=b$,  $X\Delta\{a\}\in\Phi${\rm)}.
  \end{Ax}
  Here~$\Delta$ denotes the symmetric difference of sets, $A\Delta B=(A\setminus B)\cup (B\setminus A)$.
\end{Def}
\begin{Ex}
There are three isomorphism classes of delta-matroids with the ground set of size~$1$:
\begin{equation}
D_{1,1}=\{\{1\};\{\emptyset\}\},\quad  D_{1,2}=\{\{1\};\{\emptyset,\{1\}\}\},\quad D_{1,3}=\{\{1\};\{\{1\}\}\};
\end{equation}\label{eD1}
any proper set system with a ground set of size~$1$ is a delta-matroid.

In our notation $D_{i,j}$  the first index $i$ denotes the number of elements in the ground set, while the second one is chosen ambiguously.
\end{Ex}
\begin{Def}
  A delta-matroid $(E; \Phi)$ is said to be {\em even} if for any pair of feasible sets $X$ and $Y$ {\rm(}that is, $X, Y\in\Phi${\rm)} we have
  $|X| \equiv |Y| \mod 2$. {\rm(}Here $|A|$ denotes the cardinality of a finite set $A${\rm)}.
\end{Def}

To each framed graph $G$, a delta-matroid is associated. The elements of the ground set $E$ correspond to the vertices of G,
$E=V(G)$.
 The feasible sets are best defined in terms of the adjacency matrix  $A(G)$ of a framed graph~$G$.

\begin{Def}\label{sec:--graph-to-d-matroid}
  Let $A$ be an arbitrary symmetric $|E|\times|E|$-matrix over $\mathbb{F}_2$
  whose columns and rows are marked by the elements of a finite set~$E$.\\
  Denote by~$D(A)$ the set system $D(A)=(E;\Phi(A))$ defined as follows:
  a set $F\subset E$ is feasible iff
  $\det(A|_{F})=1$ (here $A|_{F}$ is the restriction of the matrix $A$ to the subset $F$ of rows and columns).

\end{Def}

  For brevity, we write  $D(G)$ for $D(A(G))$, for a framed graph $G$.
    Bouchet proved that for every symmetric matrix $A$ with all the entries from $\FF_2$ the set system
    $D(A)$ is a delta-matroid.

Delta-matroids of the form $D(A)$ are said to be {\em graphical}.
In Eq.~(\ref{eD1}), the delta-matroids~$D_{1,1}$ and $D_{1,2}$ are graphical,
while the third one is not. The delta-matroid~$D_{1,1}$ corresponds to the
only non-framed graph on a single vertex, while the delta-matroid~$D_{1,2}$
corresponds to the one-vertex graph whose only vertex is framed.

Note that the delta-matroid $D(G)$ is even if and only if the framing of each vertex in~$G$ is~$0$
(in other words, if~$G$ is a non-framed graph).


\begin{Def}[Local duality]
  For a delta-matroid $D=(E; \Phi)$  and a subset $A\subset E$,
  define the {\em locally dual to  $D$ around $A$} by the equation
  $D*A=(E; \Phi*A)$, where
  $\Phi*A = \{F\Delta A| F\in \Phi\}$.
\end{Def}

\begin{Def}
  A delta-matroid $D=(E;\Phi)$ is said to be {\em binary} provided there is a subset $F\subset E$,
  such that $D*F$ is a graphical delta-matroid.
\end{Def}

\begin{Le}[Lemma~10~in~\cite{moffatt2017handle}]
  \label{sec:--Th-bin-with-empty-set-is-graphic}
  A binary delta-matroid is graphical if and only if the empty set is feasible.
\end{Le}

Both isomorphism classes of binary delta-matroids and even binary delta-matroids span a Hopf algebra~(\cite{LZ}).
Namely, consider the graded vector spaces
\begin{eqnarray*}
\cB&=&\cB_0\oplus\cB_1\oplus\cB_2\oplus\dots\\
\cB^e&=&\cB^e_0\oplus\cB^e_1\oplus\cB^e_2\oplus\dots,
\end{eqnarray*}
where the vector space~$\cB_n$ (respectively, $\cB^e_n$) is spanned over~$\CC$
by the isomorphism classes of binary delta-matroids (respectively,
by the isomorphism classes of even binary delta-matroids) on~$n$-element sets, $n=0,1,2,\dots$.

Multiplication in both Hopf algebras is defined as the disjoint union of set systems,
while comultiplication acts on a set system $D=(E;\Phi)$ as
$$
\mu:D\mapsto \sum_{U\sqcup W=E} D|_U\otimes D|_W.
$$
Here $D|_U$ denotes the restriction of the set system~$D$ to a subset $U\subset E$
of the ground set~$E$.

\begin{Rem}
It is shown in~\cite{Z} that the Hopf algebra~$\cB$ of binary delta-matroids is naturally
isomorphic to the Hopf algebra of Lagrangian spaces in vector spaces over~$\FF_2$ introduced in~\cite{KS}.
\end{Rem}

\subsubsection{Delta-matroids of embedded graphs}
\begin{Def}
      An {\em embedded graph\/} is  a graph drawn on a compact surface in
  such a way that its complement is a disjoint union of disks.
  \end{Def}

  We will always
  assume that the graph is connected. Edges in an embedded graph are also
  called ribbons, or handles, and we make no distinction between embedded and
  ribbon graphs. Interested reader may refer for example to~\cite{LZ} for more details on embedded graphs.\\
 Given an embedded graph $\Gamma$, one can associate to it a delta-matroid
 $D(\Gamma) = (E(\Gamma); \Phi(\Gamma))$ in a canonical way. The elements of the ground set correspond to the edges $E(\Gamma)$.
 A subset $\phi\subset E(\Gamma)$ is feasible, $\phi\in\Phi(\Gamma)$, if the boundary of the
  embedded spanning subgraph of $\Gamma$ formed by the set $\phi$ is connected, that is,
  consists of a single connected component. This means, in particular, that
  the spanning subgraph of $\Gamma$ formed by the set $\phi$ is connected (otherwise,
  each connected component would add at least one connected component to
  the boundary). Since, for a plane graph, this requirement coincides with the
  requirement that $\phi$ is a spanning tree, feasible sets for graphs embedded into
  a surface of arbitrary genus are called {\em quasi-trees}. For graphs embedded in
  surfaces of positive genus, not all the quasi-trees necessarily are trees, although
  each subset of edges forming a spanning tree is feasible.
  The pair $(E(\Gamma);\Phi(\Gamma))$, as shown by A.~Bouchet, is a binary delta-matroid.
  This binary delta-matroid is even if and only if the embedded graph~$\Gamma$ is orientable.

\subsection{Combinatorial Hopf algebras}
\label{CHA}
For definiteness, we consider Hopf algebras over the field of complex numbers~$\CC$.
A {\em combinatorial Hopf algebra\/}~\cite{ABS} is a pair $(\cH,\xi)$ consisting of a graded connected Hopf algebra~$\cH$,
$$
\cH=\cH_0\oplus\cH_1\oplus\cH_2\oplus\dots,
$$
and a linear multiplicative mapping $\xi:\cH\to\CC$.
Such mappings are called {\em characters\/} of the graded Hopf algebra~$\cH$.
 In the present paper, all the Hopf algebras
we consider are commutative and cocommutative.
For two combinatorial Hopf algebras $(\cH,\xi)$ and $(\cH',\xi')$,
a mapping $F:\cH\to\cH'$ is called a combinatorial Hopf algebras {\em morphism\/}
if it is a graded Hopf algebras morphism and takes~$\xi$ to~$\xi'$,
$\xi(F(h))=\xi'(h)$ for all $h\in\cH$, that is if the following diagram
commutes:
$$
  \xymatrix{
    {\cH}\ar[rr]^{F}\ar[dr]_{\xi} && {\cH'} \ar[dl]^{\xi'}\\
    & \CC &&
  }
$$

The ring of polynomials in infinitely many variables $\cP=\CC[p_1,p_2,\dots]$
is, in fact, a graded Hopf algebra, and we call its character $\zeta:\cP\to\CC$
defined by the requirement $\zeta:p_k\mapsto1$, $k=1,2,3,\dots$, the {\em canonical character\/}.
In the expansion
$$
\cP=\cP_0\oplus\cP_1\oplus\cP_2\oplus\dots,
$$
the space~$\cP_n$, $n=0,1,2,\dots$, is spanned by the monomials of quasihomogeneous degree~$n$;
its dimension $\dim~\cP_n$ is equal to the number of partitions of~$n$.
The comultiplication $\mu:\cP\to\cP\otimes\cP$ is the multiplicative map uniquely
determined by its values on the generators, $\mu:p_i\mapsto 1\otimes p_i+p_i\otimes1$.

The main result of~\cite{ABS}, when restricted to the case of commutative Hopf algebras,
is the following universality theorem.

\begin{Th}\label{th-ABS}
Any commmutative cocommutative combinatorial Hopf algebra~$(\cH,\xi)$ has a unique
morphism to the combinatorial Hopf algebra $(\cP,\zeta)$.
This morphism $\Psi:(\cH,\xi)\to(\cP,\zeta)$, for $h\in\cH_n$, has the form
$$
\Psi(h)=\sum_{a\vdash n}\xi^{(a)}(h)m_a,
$$
the mapping $\xi^{(a)}$ being the composition
\begin{equation}\label{e-xia}
\cH{\overset{\Delta^{k-1}}\longrightarrow}\cH^{\otimes k}\longrightarrow\cH_{a_1}\otimes\dots\otimes\cH_{a_k}{\overset{\zeta^k}\longrightarrow}\CC,
\end{equation}
where the unlabeled map is the tensor product of the canonical projections onto the homogeneous
components~$\cH_{a_i}$,
and, for $a=(a_1,\dots,a_k)$, $m_a=m_a(p_1,p_2,\dots)\in\cP$ is the
monomial symmetric function defined in the variables~$C$ as
\begin{equation}\label{e-m}
m_{a_1,\dots,a_k}\left(\sum c_i,\sum c_i^2,\sum c_i^3,\dots\right)=
\sum c_{i_1}^{a_1}c_{i_2}^{a_2}c_{i_3}^{a_3}\dots c_{i_k}^{a_k},
\end{equation}
$(i_1,i_2,\dots,i_k)$ varying along all $k=\ell(a)$-tuples of pairwise distinct integers.
\end{Th}


\section{Chromatic symmetric functions}\label{Chr}
Stanley's symmetrized chromatic polynomial of a graph can be defined in numerous ways.
We provide first the initial combinatorial definition, then the one using characters.
We use this second definition to define an extension of the symmetrized chromatic function to binary delta-matroids.
\subsection{Symmetrized chromatic polynomial of graphs}
\subsubsection{The first definition}
Let~$G$ be a simple graph. The {\em chromatic polynomial\/} of~$G$,
denoted $\chi_G(t)$, is the polynomial whose value at~$t$ is the number of proper
colorings of the vertices~$V(G)$ into~$t$ colors, $t=0,1,2,\dots$.
(A coloring is said to be {\em proper\/} if any two vertices connected by an edge
are colored in distinct colors).

Now suppose we have infinitely many colors $C=\{c_1,c_2,c_3,\dots\}$. As above, an association
$f:V(G)\to C$ to all the vertices of a simple graph~$G$ one of these colors is called a {\em proper coloring\/}
if any two vertices connected by an edge
are colored into different colors. To each coloring $f:V(G)\to C$ of the set $V(G)$ of vertices of a graph~$G$ we
associate the monomial $\gamma_f$ in the variables $c_1,c_2,c_3,\dots$, of degree $|V(G)|$,
which is equal to the product of the colors associated to the vertices,
$\gamma_f=\prod_{v\in V(G)}f(v)$.
By definition, {\em Stanley's symmetrized chromatic polynomial\/} of~$G$,
denoted $S_G(c_1,c_2,\dots)$, is the sum
$$
S_G(c_1,c_2,\dots)=\sum_{f:V(G)\to C{\text{ proper}}}\gamma_f,
$$
where the summation on the right is carried over all the proper colorings~$f$ of~$V(G)$.

Clearly, Stanley's symmetrized chromatic polynomial is a sum of infinitely many (provided~$G$
is nonempty) monomials of degree $|V(G)|$, which is symmetric under permutations of the colors.
The ring of symmetric functions of bounded degree in the variables $c_1,c_2,\dots$
can be endowed with a set of generators in a variety of ways.
In particular, it is isomorphic to the ring of
polynomials  $\CC[p_1,p_2,\dots]$ in the power sums
$$
p_k=c_1^k+c_2^k+\dots,\qquad k=1,2,3,\dots.
$$
We denote the Stanley polynomial of a graph~$G$ written in these generators by $W_G(p_1,p_2,\dots)$.
This polynomial is quasihomogeneous
of quasihomogeneous degree $|V(G)|$, if we set the degree of the variable~$p_k$
equal to~$k$ for each~$k$. Substitution $p_k=t$, for $k=1,2,3,\dots$,
makes it into the ordinary chromatic polynomial
$$
W_G(t,t,t,\dots)\equiv\chi_G(t).
$$

\begin{Rem}
A polynomial graph invariant essentially coinciding with~$W_G$ was introduced in~\cite{CDL}
for the purpose of constructing knot invariants, see below.
\end{Rem}

\subsubsection{The second definition}

Stanley's symmetrized chromatic polynomial $G\mapsto W_G(p_1,p_2,\dots)$
can be considered as a graded Hopf algebra morphism $W:\cG\to\cP$
from the Hopf algebra of graphs~$\cG$ to~$\cP$. This morphism is the one
associated, by Theorem~\ref{th-ABS}, with the character~$\xi_{\cG}:\cG\to\CC$ of~$\cG$
that takes edgeless graphs to~$1$ and graphs having edges to~$0$.
It is easy to show that the character $\xi_\cG$ indeed produces Stanley's
symmetrized chromatic polynomial. To do this, note that, for a graph~$G$
with $n=|V(G)|$ vertices, and a partition $a\vdash n$ the value
$\xi^{(a)}_{\cG}(G)$ defined by~Eq.~(\ref{e-xia}) is the number of ways to split the
set~$V(G)$ of vertices of~$G$ into disjoint parts of cardinalities $a_1,\dots,a_k$
such that the restriction of~$G$ to each of the parts is a discrete (that is, edgeless) graph.
Each such part can be colored in one of the colors $c_1,c_2,\dots$ of its own.
After summing the monomials $\gamma_f$ over all such colorings~$f$ we obtain exactly
the symmetric function~$m_a$ given by~Eq.~(\ref{e-m}).
Conversely, any proper coloring of the vertices of~$G$ determines a splitting
of the set~$V(G)$ into disjoint parts whose number is the number of distinct colors used.

\subsection{Four-term relations for graphs, weighted graphs,
and the symmetrized chromatic polynomial of framed graphs}

We are going to define the symmetrized chromatic polynomial of framed graphs
by means of a character~$\xi_{\cG^f}:\cG^f\to\CC$ that takes value~$0$ on any framed graph having
edges. Such a character is totally determined by its values
on the two framed graphs having a single vertex.
We require that the restriction of the polynomial to non-framed graphs
coincides with Stanley's symmetrized chromatic polynomial~$W$,
meaning that its value on the only non-framed graph with a
single vertex is~$1$. Our goal is to show that if the extended chromatic
polynomial satisfies the $4$-term relations for framed graphs,
as Stanley's chromatic polynomial does for non-framed ones,
then its value on the graph with a single vertex, the framing of
the vertex being~$1$, must be~$-1$.

We say that a graph invariant~$F$ {\em satisfies the $4$-term relations for
graphs}~\cite{L00} if for any graph~$G$ and any ordered pair~$(a,b)$
of its distinct vertices we have
$$
F(G)-F(G'_{ab})=F(\widetilde{G}_{ab})-F(\widetilde{G}'_{ab}).
$$
Here, all the four graphs have the same set of vertices, $V(G)$.
In the graph $G'_{ab}$ the vertices~$a,b$ are connected by an edge
if and only if they are not adjacent in~$G$ (that is, the adjacency
of~$a$ and~$b$ is switched). The graph
$\widetilde{G}_{ab}$ is obtained from~$G$ by switching the
adjacency to~$a$ of all the vertices adjacent to~$b$ and different from~$a$.
Finally, $\widetilde{G}'_{ab}=\widetilde{({G}'_{ab})}_{ab}=(\widetilde{G}_{ab})'_{ab}$.

The {\em $4$-term relations for framed graphs}~\cite{L06}
have exactly the same form, but the operation $G\mapsto \widetilde{G}_{ab}$
is understood differently if the framing of the vertex~$b$ is~$1$.
In this case, in addition to switching of the adjacency of~$a$
to the vertices adjacent to~$b$,  the framing of the vertex~$a$ is switched
as well as its adjacency to~$b$.

\begin{Le}\label{l4tch}
If~$\xi_{\cG^f}:\cG^f\to\CC$ is a character of the Hopf algebra~$\cG^f$ of framed
graphs that takes the one-vertex graph with the vertex of framing~$0$
to~$1$, takes any framed graph having edges to~$0$ and satisfies the
$4$-term relations for framed graphs, then the value of $\xi_{\cG^f}$
on the graph with a single vertex of framing~$1$ is either~$0$ or~$-1$.
\end{Le}

In order to prove the Lemma it suffices to consider the only nontrivial $4$-term relation
for framed graphs with two vertices:

\begin{picture}(400,20)(70,0)
\put(100,5){\circle*{5}} \put(140,5){\circle*{5}}
\put(200,5){\circle*{5}} \put(240,5){\circle*{5}}
\put(300,5){\circle*{5}} \put(340,5){\circle*{5}}
\put(400,5){\circle*{5}} \put(440,5){\circle*{5}}
\thicklines
\put(100,5){\line(1,0){40}} \put(400,5){\line(1,0){40}}
\put(98,11){\scriptsize{$0$}} \put(138,11){\scriptsize{$1$}}
\put(198,11){\scriptsize{$0$}} \put(238,11){\scriptsize{$1$}}
\put(298,11){\scriptsize{$1$}} \put(338,11){\scriptsize{$1$}}
\put(398,11){\scriptsize{$1$}} \put(438,11){\scriptsize{$1$}}
\put(168,4){\scriptsize{$-$}} \put(268,4){\scriptsize{$=$}}
\put(368,4){\scriptsize{$-$}}
\end{picture}

Since the value of $\xi_{\cG^f}$ vanishes in the first summand on the left
and the second summand on the right, taking into account that it is multiplicative,
we conclude the desired.

Setting the value of $\xi_{\cG^f}$ on the graph with a single vertex of framing~$1$ equal to~$0$
would lead to a character of the Hopf algebra of framed graphs that is~$0$
on any graph with at least one vertex having framing~$1$.
In order to avoid this degeneracy, we define the character
$\xi_{\cG^f}:\cG^f\to\CC$ of the Hopf algebra~$\cG^f$ as follows:
\begin{itemize}
\item the value of the character $\xi_{\cG^f}$ on the framed graph with a single vertex, whose framing is~$0$, equals~$1$;
\item the value of the character $\xi_{\cG^f}$ on the framed graph with a single vertex, whose framing is~$1$, equals~$-1$;
\item the value of the character $\xi_{\cG^f}$ on any connected framed graph with more than one vertex equals~$0$;
\item the character~$\xi_{\cG^f}$ is extended to linear combinations of graphs by multiplicativity and linearity.
\end{itemize}

Now we can define the {\em chromatic polynomial of framed graphs\/} as the morphism
$W^f:\cG^f\to\CC[p_1,p_2,\dots]$ of the combinatorial Hopf algebra~$\cG^f$ associated to the character~$\xi_{\cG^f}$.
Below, we show that the homomorphism~$W^f$ satisfies
the $4$-term relations for framed graphs.

\subsection{Chromatic polynomial for delta-matroids}
We turn  the graded Hopf algebras~$\cB$ and $\cB^e$ defined in
Sec.~\ref{sec:--simple-grp-and-even-binary-d-matroids} into combinatorial Hopf algebras by defining
the characters in the following way. Let $x\in\CC$ be a complex number.
\begin{Def}
\label{chromd}
Define $\xi_\cB:\cB\to\CC$, $\xi_{\cB^e}:\cB^e\to\CC$ by setting
\begin{itemize}
\item $\xi_{\cB}(D_{1,1})=\xi_{\cB^e}(D_{1,1})=\xi_\cB((\{1\};\{\emptyset\}))=\xi_{\cB^e}((\{1\};\{\emptyset\}))=1$;
\item $\xi_{\cB}(D_{1,2})=\xi_{\cB^e}(D_{1,2})=\xi_\cB((\{1\};\{\{1\}\})=\xi_{\cB^e}((\{1\};\{\{1\}\}))=x$;
\item $\xi_{\cB}(D_{1,3})=\xi_{\cB}((\{1\};\{\emptyset,\{1\}\}))=-1$,
\end{itemize}
and equal to~$0$ on all connected binary delta-matroids whose ground set contains more than one
element (a set system is said to be {\em connected\/} provided it is not a product
of set systems of smaller grading). As usual, the character is extended to the whole Hopf algebra
by multiplicativity and linearity.
\end{Def}


 In other words, the character $\xi_\cB$ can be computed as follows:
  $\xi(D)=0$ if $D$ is not totally disconnected;
  $\xi(D)=(-1)^kx^\ell$ if $D$ is a product of some number of~$D_{1,1}$ times $D_{1,2}^kD_{1,3}^\ell$.
(A delta-matroid $D$ is said to be {\em totally disconnected} if it can be presented as a product
  of delta-matroids having the ground set of cardinality 1).
Now we define the {\em chromatic polynomial of the Hopf algebra $\cB$} corresponding to the parameter~$x$
 as the Hopf algebra homomorphism $S_x:\cB\to\CC[p_1,p_2,\dots]$
associated to the character $\xi_{\cB}$.


Each of the three Hopf algebras $\cG,\cG^f,\cB^e$ is a Hopf subalgebra in the Hopf algebra~$\cB$.
Moreover, the following inclusions take place:
$$
  \xymatrix{
    {\cG}\ar[rr]\ar[d]\ar[drr] && {\cG^f} \ar[d]\\
     {\cB^e}\ar[rr]&& {\cB}}
  $$
The characters of these three Hopf algebras are, in fact, the restrictions of the character~$\xi_\cB$
to the corresponding subalgebras.

This definition allows one to define the chromatic polynomial of an arbitrary embedded graph:
just take for this chromatic polynomial the chromatic polynomial of the binary delta-matroid
associated to the embedded graph.
  All the three proper binary delta-matroids with a one-element ground set we have considered are,
  in fact, delta-matroids of embedded graphs. Namely, the delta-matroid $D_{1,1}=(\{1\};\{\emptyset\})$
  corresponds to the orientable embedded graph with a single vertex and one edge, the delta-matroid
  $D_{1,2}=(\{1\};\{\emptyset,\{1\}\})$ corresponds to the non-orientable embedded graph with a single vertex
  and one edge; and the delta-matroid $D_{1,3}=(\{1\};\{\{1\}\})$ corresponds to the
  embedded graph with two vertices and a single edge (note that this graph necessarily is orientable).

\section{4-term-relations for binary delta-matroids and weight systems}
\label{4T}

In~\cite{V}, V.~A.~Vassiliev introduced the notion of knot invariant of order at most~$n$,
for $n=0,1,2,\dots$. He also associated to each knot invariant of order at most~$n$
a function on chord diagrams with~$n$ chords and showed that such a function
satisfies the so-called $4$-term relations. In~\cite{Kon}, M.~Kontsevich proved
that any function on chord diagrams with~$n$ chords satisfying the $4$-term relations
arises from a knot invariant of order at most~$n$.

In our setting, a chord diagram with~$n$ chords
is nothing but an orientable embedded graph with a single vertex and~$n$ ribbons.
To each chord diagram, a simple graph is associated; this is the \emph{ intersection
graph of the chord diagram}. Its vertices are the chords of the diagram, and two
vertices are connected by an edge if the corresponding chords intersect one another.
The $4$-term relations for graphs~\cite{L00} match those for chord diagrams, and
Stanley's symmetrized chromatic polynomial satisfies these relations, whence
defining a knot invariant.

Vassiliev's $4$-term relations can be written out for arbitrary ribbon graphs,
not necessarily orientable or having a single vertex.
In the general case, they are used to describe finite order invariants of links in~$S^3$
rather than just knots; the number of vertices of a ribbon graph coincides
with the number of link components. In~\cite{LZ}, $4$-term relations were
extended to binary delta-matroids.
These $4$-term relations admit a restriction to
all the three Hopf subalgebras $\cG,\cG^f,\cB^e$ of~$\cB$, and in the case of~$\cG,\cG^f$
they coincide with the ones defined above.
Our goal is to show that the chromatic polynomial
of binary delta-matroids we have just introduced satisfies the $4$-term relations
and defines thus a link invariant.


\subsection{Vassiliev's moves for binary delta-matroids}


In order to define $4$-term relations for binary delta-matroids, we need to define two
Vassiliev's moves, the first, and the second ones.

Let~$D=(E;\Phi)$ be a binary delta-matroid, and let~$a,b\in E$ be two distinct elements of
its ground set.

\begin{Def}[The second Vassiliev move]
The result of application of the \emph{second Vassiliev move\/}
with respect to the pair $(a,b)$ to the binary delta-matroid~$D$
is defined to be
  $\widetilde D_{ab}=(E;\widetilde\Phi_{ab})$, where
  \[
    \widetilde\Phi_{ab}=\Phi\Delta\{F\sqcup\{a\}|F\sqcup\{b\}\in\Phi\text{ and } F\subset E\setminus\{a,b\}\}.
  \]
\end{Def}

The second Vassiliev move for ribbon graphs, when applied to a ribbon graph with two ribbons $a,b$
having neighboring ends, consists in the sliding of the end of the ribbon~$a$ along the ribbon~$b$.
The expression of this move in terms of the delta-matroids corresponding to the ribbon graphs
was elaborated in~\cite{moffatt2017handle}.

Alternatively, the second Vassiliev move can be defined by
\begin{Def}[Alternative definition]
Let  \begin{eqnarray*}
    \widetilde\Phi_{ab}=
                       &\{(F\not\ni a\text{ or } F\ni b)\text{ and }F\in\Phi\}\cup\\
                       &\cup\{(F\ni a\text{ \& } F\not\ni b)\text{ and }(F\in\Phi\text{ \& }F\Delta\{a, b\}\not\in\Phi)\}\cup\\
                       &\cup\{(F\ni a\text{ \& } F\not\ni b)\text{ and }(F\not\in\Phi\text{ \& }F\Delta\{a, b\}\in\Phi)\}.
  \end{eqnarray*}
\end{Def}

\begin{Def}[The first Vassiliev's move]
The result of application of the \emph{first Vassiliev move\/} to the binary delta-matroid~$D$
with respect to the pair $(a,b)$ is defined to be
  $D'_{ab}=(E;\Phi'_{ab})$, where
  \begin{eqnarray*}
   \Phi'_{ab}=&\{F\not\supset\{a, b\}\text{ and }F\in\Phi\}\cup\\
              &\cup\{F\supset\{a, b\}\text{ and }F\in\Phi\text{ and }F\setminus\{a, b\}\not\in\Phi\}\cup\\
              &\cup\{F\supset\{a, b\}\text{ and }F\not\in\Phi\text{ and }F\setminus\{a, b\}\in\Phi\}.
  \end{eqnarray*}
\end{Def}

The first Vassiliev move for ribbon graphs, when applied to a ribbon graph with two ribbons $a,b$
having neighboring ends, consists in exchanging the neighboring ends of the ribbons~$a$ and~$b$.

  For a fixed pair of elements $a$, $b$ each of the Vassiliev's moves is an involution.

    The first and the second Vassiliev's moves on the same elements pairwise commute.
    ($(\widetilde{D_{ab}})'_{ab}=\widetilde{(D'_{ab})_{ab}}$), 
which allows us to write $\widetilde{D'_{ab}}$ instead of $(\widetilde{D_{ab}})'_{ab}$.


\subsection{4-term relations for binary delta-matroids}

\begin{Def}[4-term relation for binary delta-matroids]
We say that a linear function~$\xi:\cB\to R$, where~$R$ is a commutative algebra over~$\CC$, satisfies
the \emph{$4$-term relation for binary delta-matroids\/} if
  \[
    \xi(D)-\xi(D'_{ab})=\xi(\widetilde{D_{ab}})-\xi(\widetilde{D'_{ab}}),
  \]
for arbitrary binary delta-matroid $D=(E;\Phi)$, and arbitrary pair~$a,b\in E$
of elements in its ground set.
  \label{4trel}
\end{Def}

  \begin{Th}\label{sec:4-term-relation-theorem1}
    Let $\xi:\cB\to\CC$ be a character on binary delta-matroids satisfying the 4-term relations.
    Then the Hopf algebras homomorphism $\Psi_\xi:\cB\to\cP$ also
    satisfies the 4-term relation.
\end{Th}

\begin{Proof} 
It is shown in~\cite{LZ} that the quotient~$\cF\cB$ of the Hopf algebra~$\cB$ modulo the $4$-term relations
also is a (graded commutative cocommutative) Hopf algebra. Since the character $\xi$
satisfies the $4$-term relations, its pushforward
to the Hopf algebra~$\cF\cB$ is a character of the latter algebra
and makes it into a combinatorial Hopf algebra. It defines, therefore, a Hopf algebra
homomorphism $\cF\cB\to\cP$, which is the pushforward of the homomorphism~$\Psi_\xi$,
and the theorem is proved.
\end{Proof}

Note that the restriction of a character~$\xi:\cB\to\CC$ satisfying  $4$-term relations
to any of the three Hopf subalgebras $\cG,\cG^f,\cB^e$ of~$\cB$ determines a character
of the corresponding subalgebra. This character descends to a character of the quotient of the
corresponding Hopf subalgebra modulo the $4$-term relations, and hence determines
a homomorphism of this quotient to~$\cP$.

\subsection{Chromatic weight systems}
\begin{Prop}
  \label{sec:4-term-relation-rem1}
The character $\xi_{\cB}$ defined by Def.~\ref{chromd} satisfies the 4-term relation, given by (\ref{4trel}), on binary delta-matroids:
\begin{equation}
     \xi_{\cB}(D)-\xi_{\cB}(D'_{ab})=\xi_{\cB}(\widetilde{D_{ab}})-\xi_{\cB}(\widetilde{D'_{ab}})
     \label{zeta4tr}
 \end{equation}
 \end{Prop}
  \begin{Proof}
  In the space~$\cB_1$ spanned by binary delta-matroids with a ground set of cardinality~$1$, the $4$-term relations
  are empty, hence trivially satisfied for an arbitrary linear mapping. Now suppose~$D$ is a delta-matroid
  on an $n$-element ground set, $n>1$. The value of the character~$\xi_\cB$ on a delta-matroid
 is non-zero if and only if this delta-matroid is totally disconnected,
  that is, it is a product delta-matroids on $1$-element sets. Hence, if all the four delta-matroids
  in a $4$-term relation are not totally disconnected, it is trivially satisfied by~$\xi_\cB$.
  If there is a totally disconnected delta-matroid among the four terms, then, without loss of
  generality, we can take it for~$D$.

A $4$-term relation corresponding to an ordered pair $(a,b)$ of elements of the ground
set of a totally disconnected binary delta-matroid~$D$ is nontrivial if and only if~$b$
is the element of the ground set of a factor~$D_{1,2}$ in~$D$,
while~$a$ is the element of the ground set either in a factor~$D_{1,1}$
or in a factor~$D_{1,2}$. Without loss of generality we can suppose that
the first case takes place. But the fact that the character $\xi_\cB$
satisfies this $4$-term relation has been already checked in the proof of
Lemma~\ref{l4tch}.
  \end{Proof}

\begin{Corollary}[from the theorem \ref{sec:4-term-relation-theorem1} and the proposition ~\ref{sec:4-term-relation-rem1}]\label{sec:4---1Th2}
  The extended Stanley's chromatic polynomial satisfies the 4-term relation for binary delta-matroids.
\end{Corollary}


\subsection{On the values of the extended symmetrized chromatic polynomial}
\label{Prim}
\begin{figure}
 \centering
\resizebox{.75\textwidth}{!}{\includegraphics{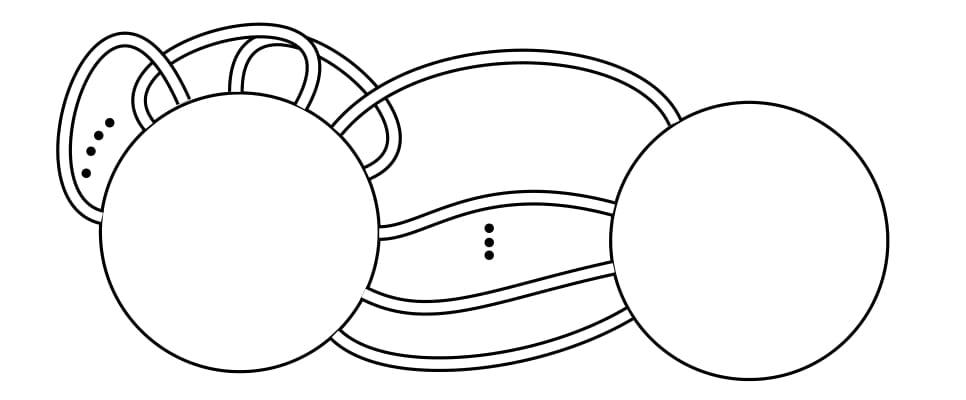}}
  \caption{A family of (orientable) embedded graphs with two vertices}\label{fff1}
\end{figure}

In the present section our goal is to estimate the distinguishing power of
the extended symmetrized chromatic polynomial~$S_x$ as a link invariant.
Being a graded Hopf algebra homomorphism,
the function~$S_x$ takes a primitive element of~$\cB_n$ to~$p_n$ times a polynomial in
the indeterminate~$x$.
We are going to show that for this polynomial we can obtain an arbitrary polynomial
of degree~$n$ whose coefficient of the linear term is~$0$.

\begin{Th}
For $n\ge2$, any polynomial in~$x$ of degree at most~$n$ with zero linear term
can be obtained as a coefficient of~$p_n$
under the mapping~$S_x$ of a primitive element in~$\cB_n$.
\end{Th}

\begin{Corollary}
The dimension of the subspace of primitive elements in~$\cB_n$ is at least~$n$.
\end{Corollary}

\begin{Proof} The coefficient of~$p_n$ in the value of~$S_x$ on a binary delta-matroid $D\in\cB_n$
coincides with the value of the character $\xi_\cB$ on the projection of~$D$
to the subspace of primitive elements along the subspace of decomposable elements.
Consider the delta-matroid of the (orientable) embedded graph with~$n$ ribbons from the family shown
in~Fig.~\ref{fff1}. If this embedded graph has~$n-k$ loops and~$k\ge2$ edges that are not loops,
then the value of~$\xi_\cB$ on its projection to the subspace of primitive elements is
proportional to~$x^k$, with a nonzero proportionality coefficient,
which proves the theorem. Note that if the number~$k$ of edges that are not loops is $k=1$,
then the delta-matroid of the embedded graph in Fig.~\ref{fff1} is a product of two
delta-matroids, one of degree $n-1$, and the other one of degree~$1$, and for $n\ge2$ its projection
to the subspace of primitive elements along the subspace of decomposable elements vanishes.

Now let us prove that the linear polynomial~$x$ cannot arise as the value of~$S_x$
on a homogeneous primitive element of degree $n\ge2$. It suffices to take into account
that the value of the character $\xi_\cB$ on an even binary delta-matroid is
non-zero only if this delta-matroid is the product of one-element binary delta-matroids.
The number of binary delta-matroids~$D_{1,2}$ among the factors in this product can be
either~$0$ or greater than~$1$, which proves the desired assertion.
\end{Proof}

\end{document}